\documentclass[11pt]{article}
\usepackage{amsmath,amssymb,amsthm}
\usepackage[mathscr]{eucal}
\usepackage{dsfont}
\usepackage{accents}
\usepackage{color,graphics,srcltx}

\newtheorem{theorem}{Theorem}[section]

\newtheorem{proposition}[theorem]{Proposition}

\newtheorem{lemma}[theorem]{Lemma}

\newtheorem{remark}[theorem]{Remark}

\def\la{{\langle}}
\def\ra{{\rangle}}
\def\IC{{\mathbb{C}}}
\def\IR{{\mathbb{R}}}

\def\cL{{\cal L}}

\def\cB{{\cal B}}
\def\cH{{\cal H}}
\def\cK{{\cal K}}

\def\cA{{\cal A}}

\def\cS{{\mathcal S}}
\def\cC{{\mathcal C}}
\def\cV{{\mathcal V}}
\def\BH{{\cB(\cH)}}
\def\bA{{\bf A}}
\def\bB{{\bf B}}
\def\bC{{\bf C}}
\def\bD{{\bf D}}
\def\bF{{\bf F}}
\def\bT{{\bf T}}

\def\bM{{\bf M}}

\def\ba{{\bf a}}

\def\bx{{\bf x}}
\def\bv{{\bf v}}

\def\bK{{\bf K}}
\def\bS{{\bf S}}
\def\b0{{\bf 0}}

\def\){{\right)}}
\def\({{\left(}}

\def\[{\left [}
\def\]{\right ]}

\def\diag{{\rm diag}}

\def\span{{\rm span}}
\def\){ \right)}
\def\({\left(}
\def\qed{\hfill\vbox{\hrule width 6 pt
\hbox{\vrule height 6 pt width 6 pt}} \medskip}

\def\conv{{\rm conv}}
\def\cl{{\bf cl}}

\begin{document}
\openup .7\jot
\title{Preservation of the joint essential matricial range}
\author{Chi-Kwong Li, Vern I. Paulsen, and  Yiu-Tung Poon}
\date{}
\maketitle

\begin{abstract}
Let
 $\bA = (A_1, \dots, A_m)$ be an $m$-tuple of
 elements of a unital C*-algebra $\cA$ and let $M_q$ denote the $q \times q$ matrices.
The  {\bf joint $q$-matricial range} $W^q(\bA)$ is the set of
$(B_1, \dots, B_m) \in \bM_q^m$ such that $B_j = \Phi(A_j)$ for some
unital completely positive linear map $\Phi: \cA \rightarrow M_q$. 
When $\cA= \cB(\cH)$, where $\cB(\cH)$ is the algebra of bounded linear operators 
on the Hilbert space $\cH$,
the {\bf joint spatial $q$-matricial range} $W^q_s(\bA)$ of $\bA$
is the set of
$(B_1, \dots, B_m) \in \bM_q^m$ such that $B_j$ is a compression of
$A_j$ to a $q$-dimensional subspace.
The {\bf joint essential spatial $q$-matricial range} is defined as
$$W_{ess}^q(\bA) =
\cap \{ \cl(W_s^q(A_1+K_1, \dots, A_m+K_m)): K_1, \dots, K_m
\hbox{ are compact operators}\},$$
where $\cl$ denotes the closure.
Suppose $\cK(\cH)$ is the set of compact operators in $\cB(\cH)$,
and  $\pi$ is the canonical surjection from $\cB(\cH)$ to the
Calkin algebra $\cB(\cH)/\cK(\cH)$. We prove that
$W_{ess}^q(\bA)$ is $C^*$-convex and equals the joint
$q$-matricial range
$W^q(\pi(\bA))$, where
$\pi(\bA) = (\pi(A_1), \dots, \pi(A_m))$.
Furthermore, for any positive integer $N$, we prove that
there are  self-adjoint compact operators $K_1, \dots, K_m$ such that
$$\cl\(W^q_s(A_1+K_1, \dots, A_m+K_m)\) = W^q_{ess}(\bA) \quad \hbox{ for all }
q \in \{1, \dots, N\}.$$ 
These results generalize those of Narcowich-Ward and Smith-Ward, obtained 
in the  $m=1$ case, and also generalize a result of M\"{u}ller obtained in case $m \ge 1$ and $q=1$. 
If $W_{ess}^1(\bA) $  is a simplex in $\IR^m$, then we prove that
there are self-adjoint compact operators $K_1, \dots, K_m$ such that
$\cl\(W^q_s(A_1+K_1, \dots, A_m+K_m)\) = W^q_{ess}(\bA)$ for all positive integers $q$.

\end{abstract}

AMS Classification. 47A12, 47A13, 47A20

Keywords. Numerical range, joint essential matricial range, 
unital completely positive linear 

\hskip .8in map,  Calkin algebra.

\section{Introduction}
Let $\cB(\cH)$ be the algebra of bounded linear operators acting on
an infinite dimensional  Hilbert space $\cH$ and let $\cK(\cH)$ denote the set of  compact operators in $\cB(\cH)$.
The {\bf numerical range} of $A \in \cB(\cH)$ is defined and denoted by
$$W(A) = \{\langle Ax, x\rangle: x \in \cH, \|x\| = 1\}.$$
It is a useful concept for studying matrices and operators.
The Toeplitz-Hausdorff Theorem states that this set is always convex
\cite{Hau,Toe}, i.e.  $tw_1+(1-t)w_2\in W(A)$ for all $w_1,\ w_2\in W(A)$
and $0\le t\le 1$. But, in general, it is not closed.

To study the joint behavior of multiple operators in $\cB(\cH)$,
researchers have considered the {\bf joint numerical range} of an
$m$-tuple $\bA = (A_1, \dots, A_m) \in \cB(\cH)^m$ defined by
$$W(\bA) = \{(\langle A_1x,x\rangle, \dots, \langle A_mx,x\rangle):
x \in \cH, \|x\| = 1\},$$
see \cite{BFL,LP0} and its references.

The
{\bf joint essential numerical range} of $\bA$ is defined as
$$W_{ess}(\bA) =
\cap \{ \cl(W(A_1+K_1, \dots, A_m+K_m)):
K_1, \dots, K_m \in \cK(\cH) \},$$
where $\cl$ denotes the closure.
These concepts were further extended to
the joint spatial $q$-matricial range, and the joint essential spatial
$q$-matricial range defined as follows. Let $\cV_q$ denote the set of operators
$X: \cK \rightarrow \cH$ such that $X^*X = I_{\cK}$
for some $q$-dimensional subpace $\cK$ of $\cH$.
The {\bf joint spatial $q$-matricial range} is
$$W^q_s(\bA) = \{(X^*A_1X, \dots, X^*A_mX):  X \in \cV_q\},$$
and the {\bf joint essential $q$-matricial range} is
$$W_{ess}^q(\bA) =
\cap \{ \cl(W^q_s(A_1+K_1, \dots, A_m+K_m)):
K_1, \dots, K_m \hbox{ are compact operators}\},$$
respectively.
Evidently, when $q = 1$, these concepts reduce to $W(\bA)$ and $W_{ess}(\bA)$.

A closely related concept is
the {\bf joint $q$-matricial range} of an
$m$-tuple of elements
$\bA = (A_1, \dots, A_m) \in \cA^m,$ where $\cA$ is a unital C*-algebra, which is defined as
$$W^q(\bA) = \{ (\Phi(A_1), \dots, \Phi(A_m) ): \Phi \hbox{ is a unital 
completely positive map from }\cA \hbox{ to } M_q \}.$$ If we let $\cS$ denote 
the operator system \cite{P0} that is the span of the identity, 
$\{ A_1,..., A_m \}$, and $\{ A_1^*,..., A_m^* \}$, then by Arveson's 
extension
theorem,
$$W^q(\bA) = \{ (\Phi(A_1), \dots, \Phi(A_m)):
\Phi \hbox{ is a unital
completely positive map from } \cS \hbox{ to } M_q\}.$$

Let $\pi: \cB(\cH)\rightarrow \cB(\cH)/\cK(\cH)$ be the
canonical map from $\cB(\cH)$ to the Calkin algebra $\cB(\cH)/\cK(\cH)$
and set $\pi(\bA) = (\pi(A_1), \dots, \pi(A_m))$.
Then it is clear that
$$\cl(W^q_s(\bA)) \subseteq W^q(\bA) \qquad \hbox{ and } \qquad
 W^q(\pi(\bA)) \subseteq W^q(\bA+\bK) \quad   \mbox{ for all }\quad \bK\in \cK(\cH)^m $$
In this paper, we will show that $W_{ess}^q(\bA) = W^q(\pi(\bA))$
and   consequently,
$W_{ess}^q(\bA)$ is $C^*$-convex.
 
When $m=1$, this reduces to a result of Narcowich and Ward  \cite{NW}.
Moreover, we study the preservation problem for $W_{ess}^q(\bA)$ and
$W^q(\pi(\bA))$, namely, we prove that for each $N$, there is an $m$-tuple of self-adjoint compact operator
$\bK = (K_1, \dots, K_m)$ such that for $1 \le q \le N$,
\begin{equation} \label{eq1}
\cl(W^q_s(\bA+\bK)) = W_{ess}^q(\bA)=
W^q( \bA + \bK) = W^q(\pi(\bA)).
\end{equation}
When $m=1$ this reduces to a result of Smith and Ward \cite{SW};
when $N = 1$, this reduces to a result of M\"{u}ller \cite{M}.

Let $\cS(\cH)=\{A\in \cB(\cH)):A=A^*\}$.   Note that every $A \in \cB(\cH)$
has a unique Cartesian decomposition
$A = A_1 + iA_2$ such that $A_1, A_2\in \cS(\cH)$.
Thus, one can identify $W^q_s(A)$ with $W^q_s(A_1,A_2)$, and also
identify $W_{ess}^q(A)$ with  $W_{ess}^q(A_1, A_2)$.
For this reason, we shall henceforth focus on  $W^q_s(\bA)$ and $W_{ess}^q(\bA)$ for
$\bA\in  \cS(\cH)^m$.

We show that if $W_{ess}(\bA)$ is a simplex in $\IR^m$, then
there are compact operators $K_1, \dots, K_m$ such that (\ref{eq1}) holds
for all positive integers $q$. This extends another result in \cite{SW}.

In our discussion, we will always assume that $\cH$ is infinite dimensional.
In addition to the notation $\cB(\cH), \cK(\cH)$ and $\cS(\cH)$ introduced 
above,  we will let $\cS_{\cK}(\cH)$ be the set of compact operators in 
$\cS(\cH)$.

\section{Some basic results}
\setcounter{equation}{0}

The following result extends \cite[Theorem 2.1]{Farenick2}.

\begin{theorem} \label{thm2.1} Let $\cA$ be a unital C*-algebra and 
let $\bA = (A_1, \dots, A_m)\in \cA^m$ be self-adjoint.
Then $(B_1, \dots, B_m) \in W^q(\bA)$ if and only if
\begin{equation} \label{rb}
\|R_0 \otimes I_q + R_1 \otimes B_1 + \cdots + R_m \otimes B_m\|
\le \|R_0 \otimes I + R_1 \otimes A_1 + \cdots + R_m \otimes A_m\|
\end{equation}
for all $R_0, \dots, R_m \in M_q$.
\end{theorem}

\it Proof. \rm Note that (\ref{rb}) is equivalent to the condition that 
the unital  linear map sending $A_j\in \cA$
to the Hermitian matrix $B_j\in M_q$ is completely contractive. Then the result 
follows from the fact that
every unital completely contractive map
from an operator system to $M_q$
is a unital completely positive linear map; for example,
see \cite{P0}.
\qed

Given $\bB= (B_1,...,B_m) \in M_q^m$ and $L \in M_q$ we set
$$L^* \bB L= (L^*B_1L,..., L^*B_m L).$$
A subset  $\cC \subseteq M_q^m$ is {\bf $C^*$-convex} if
$$\sum_{j=1}^N L_j ^*\bB_j L_j \in \cC$$
for any $\bB_1, \dots, \bB_N \in \cC$ and $L_1, \dots, L_N \in M_q$
satisfying $\sum_{j=1}^N L_j^*L_j = I_q$. It is well known that
$W^q(\bA)$ is $C^*$-convex. We have the following result showing that
$W_{ess}^q(\bA) = W^q( \pi(\bA))$. 
It will then follow that $W^q_{ess}(\bA)$ is also $C^*$-convex.

\begin{theorem} \label{thm2.2}
Let $\bA = (A_1, \dots, A_m) \in \cS(\cH)^m$.
Then

\begin{equation}
\label{setequality}
W_{ess}^q(\bA) = W^q(\pi(\bA)),
\end{equation}
which consists of  $\bB = (B_1, \dots, B_m) \in M_q^m$ satisfying:
\begin{equation}
\label{rba}\|R_0 \otimes I_q + R_1 \otimes B_1 + \cdots + R_m \otimes B_m\|
\le \|R_0 \otimes I + R_1 \otimes \pi(A_1) + \cdots + R_m \otimes \pi(A_m)
\|\end{equation}
for all $R_0, \dots, R_m \in M_q$. Consequently, $W_{ess}^q(\bA)$ is $C^*$-
convex. 
\end{theorem}

\medskip
\it Proof.\   \rm  Let $S$ be the finite dimensional operator system  given 
by the span of $\{\pi(I),\pi(A_1),\dots,\pi(A_m)\}$. Then by 
\cite[Theorem 9.11]{Kavruk}, there exists a unital $q$-positive map 
$R:S\to \BH$ such that $R( \pi(A_j) )=A_j+K_j$ with  
$K_j\in \cS_{\cK}(\cH)$ for $1\le j\le m$.

Let  $\bB = (B_1, \dots, B_m) \in W_{ess}^q(\bA)$. For $\varepsilon >0$,
there is
$$(\tilde{B}_1, \dots,  \tilde{B}_m)  \in W^q_s(A_1+K_1,\dots,A_m+K_m)$$
such that $\|\tilde{B}_j-B_j\|<\varepsilon $ for all $1\le j\le m$.
Let $X\in\cV_q$ be such that $X^*(A_j+F_j)X=\tilde{B}_j$ 
for all $1\le j \le m$.
Then the map $\Phi:S\to M_q$ given by $\Phi( T)=X^*R(T)X$ is a unital
completely positive map   satisfying $(\tilde{B}_1, \dots,  \tilde{B}_m)
=\(\Phi(\pi(A_1)),\dots,\Phi(\pi(A_m))\)\in  W^q(\pi(\bA))$. 
Since $ W^q(\pi(\bA))$ is close, $\bB\in W^q(\pi(\bA))$.
Hence, $W_{ess}^q(\bA) \subseteq W^q(\pi(\bA))$.

For the reverse inclusion, if $\bB=(B_1, \dots, B_m) \in W^q(\pi(\bA))$.
Let $ \phi:S\to M_q$ be a unital completely positive map 
such that $\phi(\pi(A_j))=B_j$ for
all $1\le j\le m$. Let $\Phi=\phi\circ \pi$. Then $\Phi$ is a unital 
completely positive on the separable $C$*-algebra 
$\mathcal{A}$ generated by 
$\{ I, A_1,...,A_m \}$, with $\Phi( \mathcal{A} \cap \mathcal{K}
(\mathcal{H})) = (0).$
Hence, by \cite[Theorem~2.5]{BS} (see also \cite[Theorem 4]{Ar}), we have 
that $\bB=(B_1, \dots, B_m) \in \cl(W^q_s(A_1,..., A_m)).$
Given $K_1,\dots, K_m \in\cS_{\cK}(\cH)$, we also have that 
$\Phi$ vanishes on the intersection of the compact operators with the 
separable C*-algebra generated by $\{ I, A_1+K_1, \dots A_m+K_m \}$ 
with $\Phi(A_j+K_j)=\phi(\pi(A_j)) $.
Hence, by the same reasoning,  
$\bB \in \cl(W^q_s(A_1+K_1, \dots, A_m+K_m))$.
 
This proves   (\ref{setequality}). (\ref{rba}) follows from 
Theorem \ref{thm2.1} and the last statement is a consequence of the 
$C^*$-convexity of $W^q(\pi(\bA))$.\qed
 
\begin{theorem} \label{thm2.5}
Let $\bA = (A_1, \dots, A_m)$ be an $m$-tuple of self-adjoint operators.
The following conditions are equivalent.
\begin{itemize}
\item[{\rm (a)}]
$W_{ess} (\bA) = W (\pi(\bA))$ has non-empty interior in $\IR^{1\times m}$.
\item[{\rm (b)}] For any positive integer $q$,
$W_{ess}^q(\bA) = W^q(\pi(\bA))$ has non-empty interior in $\cH_q^m$,
where $\cH_q$ is the real linear space of $q\times q$ Hermitian matrices.
\item[{\rm (c)}] The set $\{I, \pi(A_1), \dots, \pi(A_m)\}$
is linearly independent.
\end{itemize}
\end{theorem}

\it Proof. \rm (a) $\Rightarrow$ (b):
Suppose $W_{ess}(\bA)$ has non-empty interior. We may assume that 
for some $r>0$,
$D_r=\{\ba\in \IR^{1\times m}:
\|\ba\|_{\infty}\le r\}\subseteq W_{ess}(\bA)$. Let $q>1$. 
We are going to  show that
$(B_1,\dots,B_m)\in W_{ess}^q(\bA)$ for all $B_1,\dots,B_m\in H_q$ with
$\|B_i\|\le \dfrac{r}{m}$,  $1\le i\le m$.

Let $|b|\le r$. Then $(b,0,\dots,0)\in W_{ess}(\bA)$. 
So there exists \cite[Theorem 2.1]{LP1} an
orthonormal  sequence of vectors $\{\bx_n\}$ such that
$\underset{n\to \infty} \lim \(\la A_1 \bx_n,\bx_n\ra,\dots,\la
A_m\bx_n,\bx_n\ra\) = (b,0,\dots,0)$.
For every  $K\in \cK(\cH)$, we have
$\underset{n\to \infty} \lim \la  K\bx_n,\bx_n\ra= 0$.
Let $S=$ span $\{I,A_1,\dots,A_m\}$. Using
the    quotient map $\pi$ from $\cB(\cH)$ to the Calkin
algebra $\cB(\cH)/\cK(\cH)$, we have
$\phi(\pi(A))=\underset{n\to \infty} \lim \la  A\bx_n,\bx_n\ra$
defines a unital completely positive map on $\pi(S)$   such that 
$ \phi(\pi(A_1))=b$
and
$\phi(\pi(A_2))= \cdots=\phi(\pi(A_m))=0$. $\phi$ can be extended to a unital
completely positive map
on $\cB(\cH)/\cK(\cH)$.
Let  $B=U^*\diag(b_1,\dots,b_q)U\in H_q$ with $(b_1,\dots,b_q)\in D_r $. Then for each
$1\le i\le q$,
there is a unital  completely positive map $\phi_i$ on $\cB(\cH)/\cK(\cH)$  such that
$\phi_i(\pi(A_1))=b_i$ and $\phi_i(\pi(A_2))= \cdots=\phi_i(\pi(A_m))=0$. Let
$\Phi(\pi(A))=U^*\(\diag( \phi_1(\pi(A)),\dots,\phi_q(\pi(A)))\)U$.
Then $\Phi: \cB(\cH)/\cK(\cH)\to M_q$ is a
unital completely positive map such that $\Phi\(\pi(A_1)\)=B$ and
$\Phi\(\pi(A_2)\)= \cdots=\Phi\(\pi(A_m)\)=0$. Suppose $(B_1,\dots,B_m)\in H_q^m$
with $\|B_i\|\le \dfrac{r}{m}$ for all $1\le i\le m$. From the above discussion,
there exist unital completely positive maps
$\Phi_i:\cB(\cH)/\cK(\cH)\to M_q$ such that
$\Phi_i(\pi(A_j))=\delta_{i\,j}\(mB_j\)$ for
all $1\le i,\,j\le m$. Let $\Phi=\dfrac{1}{m}\(\sum_{i=1}^m\Phi_i\)$. Then $\Phi$ is a unital completely positive map satisfying
$\Phi(\pi(A_j))=B_j$ for  all $1\le j\le m$.

(b) $\Rightarrow$ (c): Suppose (c) is not true. Then there are
real numbers $a_0, \dots, a_m$ not all zero
$a_0 I + a_1 \pi(A_1) + \cdots + a_m \pi(A_m) = 0$.
Then $a_1B+ \cdots + a_m B =- a_0 I_q$ for
every $(B_1, \dots, B_q) \in W^q(\pi(\bA))$. So, $W^q(\pi(\bA))$
has empty interior.

(c) $\Rightarrow$ (a): Suppose  (a) is not ture. Then the convex set
$W( \pi(\bA))$ has empty interior. So, it must lie in an affine space,
say, $\{ (x_1, \dots, x_m) \in \IR^{1\times m}:
a_0 + a_1x_1 + \cdots + a_m x_m= 0\},$
where $(a_0, a_1, \dots, a_m)\in \IR^{1\times (m+1)}$ is a unit vector.
It then follows that $\langle\( a_0 I + \cdots a_m \pi(A_m)\)v,v\rangle=0$
for any unit vector $v$. Thus,
$a_0 I + \cdots + a_m \pi(A_m) = 0$, i.e.,
$\{I, \pi(A_1), \dots, \pi(A_m)\}$ is linearly dependent.

\qed

\section{Preservation problems}
\setcounter{equation}{0}

Narcowich and Ward \cite{NW} proved that given a single
operator $A \in \cB(\cH)$ one has $W^q(\pi(A)) = W_{ess}^q(A)$
for every $q$.
Smith and Ward \cite[Section 5]{SW} proved that  given a single
operator $A \in \cB(\cH)$
and a positive integer $N$,
there exists $K \in \cK(\cH)$ such that $W^q(A+K) = W^q( \pi(A)),$
for all $q=1,\dots,N.$ M\"{u}ller \cite[Corollary 14]{M} proves that given an  $m$-tuple of  operators $\bA = (A_1, \dots, A_m)$,  
there is an $m$-tuple of compact  operators
$\bK = (K_1, \dots, K_m)$
such that $\cl(W(\bA + \bK)) = W_{ess}(\bA)$. 
 The following result extends these results  to joint (spatial, essential) $q$-matricial ranges of tuples of operators. 
Again, we can focus on tuples of 
self-adjoint operators.

\begin{theorem} \label{thm3.1}
Suppose $\bA = (A_1, \dots, A_m)$ is an $m$-tuple of self-adjoint operators.
Then for any positive integer $N$,
there is an $m$-tuple of compact self-adjoint operators
$\bK = (K_1, \dots, K_m)$
such that
\begin{equation}\label{fourset}
\cl(W^q_s(\bA + \bK)) = W^q( \bA + \bK) = W^q(\pi(\bA)) = W_{ess}^q(\bA)
\quad \hbox{ for all } q = 1, \dots, N.
\end{equation}
\end{theorem}

\it Proof. \rm
First, we show that there is $\bK$ such that
$$W^N(\bA + \bK) = W^N(\pi(\bA)).$$
Given $\bA = (A_1,..., A_m)$,
let $S$ be the operator system of dimension
not larger than $m+1$ spanned by
$\pi(I), \pi(A_1), \dots, \pi(A_m)$.
This operator system is embedded into the Calkin algebra by the identity map.

By   \cite[Theorem 9.11]{Kavruk},
there is an $N$-positive  lifting map
$R: S\rightarrow \cB(\cH)$.  Since it is an $N$-positive     lifting,
$R(\pi(A_i)) = A_i + K_i$  for some $(K_1, \dots, K_m)\in\cS_{\cK}(\cH)$.

Since  $A_i+K_i\mapsto \pi(A_i)$ is completely positive,  
it follows that 

$$
W^q( A_1+K_1,..., A_m+K_m) = W^q( \pi(A_1),...,\pi(A_m))
\quad \hbox{ for all } q \le N.$$
To see this, let $(T_1,...,T_m)$ be in the left hand side.
Then there is a completely positive map $\phi: S_0 \rightarrow M_q$ with
$T_i=\phi(A_i+K_i) =\phi(R(\pi(A_i))) $. But then the map
$\gamma(\pi(A_i)) = T_i$ is $N$-positive.  Since $q \le N$, by \cite[Theorem~6.1]{P0},
$\gamma$ is completely positive. Hence
$(T_1,...,T_m)$ is in the right hand side.
The proof of other containment is similar.

Now, using the notation $\bA = (A_1, \dots, A_m), \bK= (K_1,\dots, K_m)$
we note that
$$W_{ess}^q(\bA) \subseteq \cl(W^q_s(\bA + \bK)) \subseteq W^q(\bA+\bK)$$
is always true. By Theorem \ref{thm2.2} and the above discussion,
$$W^q(\bA+\bK) = W^q( \pi(\bA)) = W_{ess}^q(\bA),$$
so that (\ref{fourset}) is true.
\qed

\begin{remark} \rm
The still open Smith-Ward problem asks: given $A \in \cB(\cH)$
does there exist a compact operator $K$ such that  for every $q \in \mathbb{N}$, $W^q(A+K) = W^q( \pi(A))$? Let $A=A_1+iA_2$ be the Cartesian decomposition of $A$.
It is known \cite{P} that this problem is equivalent to asking if the operator system $\cS_{\pi} := \span\{\pi(I), \pi(A_1), \pi(A_2)\}\subseteq \cB(\cH)/\cK(\cH)$
has a completely positive lifting   
to $\cB(\cH)$, i.e., if there exists a unital completely positive map $\Phi: \cS_{\pi} \to \cB(\cH)$ such that $\pi(\Phi(x)) = x$ for all $x \in \cS_{\pi}$. The analogue of this problem for tuples of operators is known
to be false.  In \cite[Theorem 3.3]{P} an example is given of an operator $A \in \cB(\cH)$
such that the span of  
$\{ \pi(I), \pi(A),  \pi(A^*),  \pi(A^*A),  \pi(AA^*) \}$ 
does not have  a completely positive lifting  to $\cB(\cH)$.
This implies that for the tuple $\bA =( A_1, A_2, A^*A, AA^*)$, 
there does not exist four compact operators $\bK =(K_1, K_2, K_3, K_4)$ such that
$W^q( \bA + \bK) = W^q( \pi(\bA))$ for every $q$. This example shows
that the finite range on $q$ in the above theorem is necessary. It also
shows that the analogue of the Smith-Ward problem is false for
four or more self-adjoint operators.
\end{remark}

In \cite[Theorem 5.1]{SW}, the authors
showed that if  $A \in \cB(\cH)$ is such that
$W_{ess} (A)$ is a line segment, then there is a compact operator
$K$ such that $W^q_s(A+K) = W_{ess}^q(A)$ for all $q = 1, 2, \dots$.
We can extend the result to the situation when $W_{ess}(\bA) =
W( \pi(\bA))$ is a simplex in $\IR^{1\times m}$. 
To achieve this, we need the following
result from \cite[Theorem 1.1]{BFL} and a lemma.

\begin{theorem} \label{thm3.2}
Let $\bT = (T_1, \dots, T_m) \in \cS(\cH)^m$, and
$\bD = (D_1, \dots, D_m) \in M_{m+1}$
be an $m$ tuple of diagonal matrices with nonnegative
diagonal entries such that
$W(\bD)$ is a simplex in $\IR^{1\times m}$.
Then  $W(\bT) \subseteq W(\bD)$ if and only if there is a Hilbert space
$\tilde \cH$ and  $X: \cH \rightarrow \tilde \cH \otimes \IC^{m+1}$
satisfying $X^*X = I_\cH$ and $T_j = X^*(I_{\tilde \cH} \otimes D_j)X$
for $j = 1, \dots, m$.
\end{theorem}

\begin{lemma} \label{lem3.4}
Let $\bA = (A_1, \dots, A_m) \in \cS(\cH)^m$. Suppose  $\bB_i=(B_{i1},\dots,B_{im})\in W_{ess}^p(\bA)$ for $1\le i\le N$. Then for all $\varepsilon>0$ and $K_1,\dots,K_m\in \cS_{\cK}(\cH)$, there exists $\tilde{\bB}=(\tilde{B}_{1},\dots,\tilde{B}_{m})\in W^{Np}(A_1+K_1,\dots,A_m+K_m)$  such that $\tilde{B}_{j}=\oplus_{i=1}^N\tilde{B}_{ij}$ with $\|B_{ij}-\tilde{B}_{ij}\|<\varepsilon$ for all  $1\le i\le N$ and $1\le j\le m$.
\end{lemma}

\it Proof. \rm  Without loss of generality, we may assume that $K_1=\cdots   =K_m=0$.
We are going to prove by induction on $N$. For $N=1$, the result follows from definition.

Suppose the result holds for some $N=k\ge 1$. Then there is
$X_1:\cK_1 \to \cH$ with $X_1^*X_1 = I_{\cK_1}$
for some $kp$-dimensional subspace $\cK_1$ of $\cH$
such that $X_1^*A_j X_1 =\oplus_{i=1}^k\tilde{B}_{ij}$ for $j =1 ,\dots,m$ and $\|B_{ij}-\tilde{B}_{ij}\|<\varepsilon$ for all  $1\le i\le k$ and $1\le j\le m$.

Extend the operator $X_1$ to an unitary operator $U:\cH \to \cH$ so that $U|_{\cK_1} = X_1$.
Let $\cL$ be the subspace spanned by
$$\left\{\cK_1,\,  U^*A_1U(\cK_1),\, \dots,\, U^*A_mU(\cK_1) \right\}.$$
Then $\cL$ has dimension at most $kp(m+1)$. Define $Y = U|_{\cL^\bot}$.
Then $Y^*Y = I_{\cL^\bot}$.
Furthermore, as $X_1  = U|_{\cK_1}$, $Y = U|_{\cL^\bot}$
and $\cK_1  \subseteq \cL$, we have $Y^*X_1 = 0$.
Also for any $u \in \cL^\bot$ and $v \in \cK_1$,
$U^*A_jUv \in \cL$. Hence,
$$\langle u, Y^*A_j X_1 v\rangle
= \langle  Yu, A_j X_1 v \rangle
= \langle  Uu, A_j U v \rangle
= \langle  u, U^*A_j U v \rangle
= 0$$
and thus, $Y^*A_j X_1 = 0$ for all $j =1,\dots,m$.

Since  $Y^*A_jY-A_j \in  \cS_{\cK}(\cH)$ for all  $1\le j\le m$, there exists
$X_2: \cK_2 \to  \cL^\bot$ with $X_2^*X_2 = I_{\cK_2}$
for some $p$-dimensional subspace $\cK_2$ of $\cL^\bot$
such that  $X_2^*(Y^* A_j Y) X_2 = \tilde{B}_{(k+1)j}$ satisfies
$$\|B_{(k+1)j}-\tilde{B}_{(k+1)j}\|<\varepsilon\mbox{ for all   }1\le j\le m.$$ 

Observe that  $\cK_1$ and $\cK_2$ are two mutually orthogonal subspaces of $\cH$.
Furthermore, $X_1: \cK_1 \to \cH$ and $YX_2: \cK_2 \to \cH$ satisfy
$$X_1^*X_1 = I_{\cK_1},\quad (YX_2)^*(YX_2) = I_{\cK_2},\quad
(YX_2)^*X_1 = 0
\quad\hbox{and}\quad
(YX_2)^*A_j X_1 = 0$$
for $ j =1,\dots,m$.
Then the operator $Z = X_1 \oplus (YX_2): \cK_1 \oplus \cK_2 \to \cH$
satisfies $Z^*Z = I_{\cK_1\oplus \cK_2}$ and
 $Z^*A_j Z =   \oplus_{i=1}^{k+1}\tilde{B}_{ij}$ for $j =1 ,\dots,m$
\qed

Now we are ready to present the other main result for this section.

\begin{theorem} \label{thm3.3} Suppose $\bA = (A_1, \dots, A_m)$
is an $m$-tuple of self-adjoint operators acting on $\cH$. If
$W_{ess}(\bA) $ is a simplex $\bS$ in $\IR^{1\times m}$, i.e.,
$W_{ess}(\bA) $ is a
polyhedron with $m+1$ vertices,
then there is an $m$-tuple of self-adjoint compact operators
$\bK = (K_1, \dots, K_m)$ such that for all positive integer $q$,
\begin{equation}
\label{fourset-2}
\cl(W^q_s(\bA + \bK)) = W^q(\bA+\bK) =  W_{ess}^q(\bA) = W^q( \pi(\bA)).
\end{equation}
\end{theorem}

\it Proof. \rm
By Theorem \ref{thm3.1}, there is an $m$-tuple of
compact operator $\bK$ such that
$$\cl(W(\bA + \bK)) = W(\bA + \bK) = W(\pi(\bA)) = W_{ess}(\bA).$$
We will show that for any positive integer $q$, we have
$$\cl(W^q_s(\bA + \bK)) = W^q(\bA + \bK) = W^q(\pi(\bA)) = W_{ess}^q(\bA).$$

Let $\bD = (D_1, \dots, D_m)$ such that
$D_i = \diag(v_{i1}, \dots, v_{i,m+1})$ for $i = 1, \dots, m$,
where 
$$\bv_1=(v_{11}, \dots, v_{m1}),\dots, \bv_{m+1}=(v_{1\,m+1}, 
\dots, v_{m\,m+1}) \in \IR^{1\times m}$$
are the vertices of the simplex $\bS$.
Then $\bS = W(D_1, \dots, D_m) = W_{ess}(\bA)$.

Suppose $q>1$ is a positive integer.
Let
$\bB = (B_1, \dots, B_m) \in W^q(\bA+\bK)$.
Then for any unital completely positive map $\phi: M_k\rightarrow M_1$,
$(\phi(B_1), \dots, \phi(B_m)) \in  \cl(W(\bA + \bK))$, which equals to
the simplex $\bS$. Thus, $W(B_1, \dots, B_m) \subseteq W(D_1, \dots, D_m)$.
By Theorem  \ref{thm3.2}, $(B_1, \dots, B_m)$ admits a joint dilation
to $(I\otimes  D_1, \dots, I \otimes  D_m)$. Hence,  $W^q(B_1, \dots, B_m) \subseteq W^q(D_1, \dots, D_m)$.

Suppose $\bK = (K_1, \dots, K_m) \in \cS(\cH)^m$. 
Note that $\bv_1,\dots,\bv_{m+1}\in  W_{ess}(A_1+K_1,\dots,A_m+K_m)$.  
Applying Lemma \ref{lem3.4} 
with $p=1$, $N=(m+1)q$ and $\bB_i=\bv_j$  if $i\equiv j \ ({\rm mod\ }m+1)$ for 
$1\le i\le N$ and $1\le j\le m+1$, we can get a 
sequence $\bD^{(k)}=\(D^{(k)}_1,\dots D^{(k)}_{m}\)\in 
W^{N}_s(A_1+K_1,\dots,A_m+K_m)$ such that 
$$\{D_j^{(1)}, D_j^{(2)}, D_j^{(3)}, \dots, \} \rightarrow I_q\otimes D_j.$$

Since $W^q(B_1, \dots, B_k) \subseteq W^q(D_1, \dots, D_m)$,
for any $R_0, \dots, R_m \in M_q$ we have
 \begin{eqnarray*}&&\|R_0 \otimes I_q + R_1 \otimes B_1 + \cdots + R_m \otimes B_m\|\\
&\le& \|R_0 \otimes I_{m+1} + R_1 \otimes D_1 + \cdots + R_m \otimes D_m\|\\
&=& \|R_0 \otimes I_{m+1} \otimes I_q \otimes I_{m+1} + R_1 \otimes I_{m+1} \otimes I_q\otimes D_1 + \cdots + R_m \otimes I_{m+1} \otimes I_q  \otimes D_m\|\\
&=& \lim_{k\rightarrow \infty}
\|R_0 \otimes I_{m+1}\otimes I_{N} + R_1\otimes I_{m+1} \otimes D_1^{(k)} + \cdots + R_m\otimes I_{m+1} \otimes D_m^{(k)}\|\\
&\le&\| R_0\otimes I_{m+1} \otimes I + R_1\otimes I_{m+1} \otimes (A_1+K_1) + \cdots + R_m\otimes I_{m+1} \otimes (A_m+K_m)\|\\
&=&
\| R_0 \otimes I + R_1 \otimes (A_1+K_1) + \cdots + R_m \otimes (A_m+K_m)\|.
\end{eqnarray*}
 
By Theorem \ref{thm2.1}, $\bB \in W^q(\bA+\bK)$.
Because this is true for all compact $\bK$, we see that
$\bB \in W_{ess}^q(\bA).$
Hence, we have
$$\cl(W^q_s(\bA+\bK)) \subseteq W^q(\bA+\bK) \subseteq W_{ess}^q(\bA).$$
Since $W^q(\pi(\bA)) = W_{ess}^q(\bA+\bK)) \subseteq \cl(W^q_s(\bA+\bK))$,
we see that (\ref{fourset-2}) holds. \qed

\section{Related results}

By Theorem \ref{thm3.3}, we have
the following   an  extension of \cite[1.22.5]{Sakai}.

\begin{proposition} \label{3.5}
Let $\bA = (A_1, \dots, A_m) \in \cS(\cH)^m$. Suppose
$W_{ess}(\bA)$ is a subset of a simplex $\bS$ in $\IR^{1\times m}$ with vertices
$v_1, \dots, v_{m+1}$
such that $v_k = (v_{1k}, \dots, v_{m,k})$ for $k = 1, \dots, m+1$.
Then there is an $m$-tuple of self-adjoint compact operators
$\bK = (K_1, \dots, K_m)$
such that for any $R_0, \dots, R_m \in M_q$,
$$\ \qquad \|R_0 \otimes I + R_1 \otimes (A_1+K_1) + \cdots + R_m \otimes (A_m+K_m)\|
$$
\begin{equation} \label{norminequality}
\le \max\{ \|R_0 + v_{1k}R_1 + \cdots + v_{mk}R_m\|: 1 \le k \le m+1\}.
\end{equation}
In fact, $\bK$ can be chosen such that 
the equality holds in {\rm (\ref{norminequality})} for any choice of 
$R_0, \dots, R_m \in M_q$.
\end{proposition}

\it Proof. \rm 
Let $D_1, \dots, D_m \in M_{m+1}$ with 
$D_j = \diag(v_{j1}, \dots, v_{j,m+1})$ so that 
$W(D_1, \dots, D_m) = \bS$.
Let $\tilde \cH = \cH \oplus (\cH \otimes \IC^{m+1})$. Then 
$A_j$ is a compression of $\tilde A_j = A_j \oplus (I_\cH \otimes D_j)
\in \cB(\tilde \cH)$
for all $j = 1, \dots, m$. Evidently, 
$W_{ess}^q(\tilde A_1, \dots, \tilde A_m) = W^q(D_1, \dots, D_m) $.
By Theorem \ref{thm3.3}, there are self-adjoint  compact operators
$\tilde K_1, \dots, \tilde K_m \in \cB(\tilde \cH)$
such that $W(\tilde A_1 + \tilde K_1, \dots, \tilde A_m + \tilde K_m)
= W_{ess}(\tilde A_1, \dots, \tilde A_m)$.
Moreover, for any positive integer $q$, we have
$$W^q(\tilde A_1 + \tilde K_1, \dots, \tilde A_m + \tilde K_m)
= W_{ess}^q(\tilde A_1  , \dots, \tilde A_m  ).$$
 
Suppose $\tilde K_j$ has operator matrix 
$\begin{pmatrix} K_1 & * \cr * & * \cr\end{pmatrix}$ with $K_1 \in \cB(\cH)$
for $j = 1, \dots, m$.
Then for any $R_0, \dots, R_m \in M_q$, we have
\begin{eqnarray*}
&&\|R_0 \otimes I + R_1 \otimes (A_1+K_1) + \cdots + R_m \otimes (A_m + K_m)\| \\
&\le& \|R_0 \otimes I + R_1 \otimes (\tilde A_1+\tilde K_1) 
+ \cdots + R_m \otimes (\tilde A_m + \tilde K_m)\| \\
& = &  \|R_0 \otimes I + R_1 \otimes D_1 + \cdots + R_m \otimes D_m\|\\
&=&\max\{ \|R_0 + v_{1k}R_1 + \cdots + v_{mk}R_m\|: 1 \le k \le m+1\}.
\end{eqnarray*}

For the last assertion, let $\bK$ be chosen such that (\ref{norminequality})
holds.
Suppose $A_j + K_j$ has operator matrix
$\begin{pmatrix}B_j & * \cr * & C_j \cr \end{pmatrix}$ with $B_j \in M_{m+1}$
for $j = 1, \dots, m$. There are finite rank operators $F_1, \dots, F_m$ such that
$A_j+K_j + F_j$ has operator matrix $D_j \oplus C_j$ for $j = 1, \dots, m$.
Then for any $R_0, \dots, R_m \in M_q$, we have
\begin{eqnarray*}&&\|R_0 \otimes I + R_1 \otimes C_1 + \cdots + R_m \otimes C_m\|\\ 
&\le&
\|R_0 \otimes I + R_1 \otimes (A_1+K_1) + \cdots + R_m \otimes (A_m + K_m)\|
\\ 
&\le&\|R_0 \otimes I + R_1 \otimes D_1 + \cdots + R_m \otimes D_m\|.\end{eqnarray*}
Hence,  if $\hat K_j = K_j + F_j$ for $j = 1, \dots, m$, then
\begin{eqnarray*}
&&\|R_0 \otimes I + R_1 \otimes (A_1+\hat K_1) + \cdots + R_m \otimes 
(A_m +\hat K_m)\| \\
&=& \max\{\|R_0 \otimes I + R_1 \otimes D_1 + \cdots + R_m \otimes 
D_m\|,
\|R_0 \otimes I + R_1 \otimes C_1 + \cdots + R_m \otimes C_m\|\} \\
&=& \|R_0 \otimes I + R_1 \otimes D_1 + \cdots + R_m \otimes D_m\| \\
&=& \max\{ \|R_0 + v_{1k}R_1 + \cdots + v_{mk}R_m\|: 1 \le k \le m+1\}.
\end{eqnarray*}
The last assertion follows. \qed

Motivated by quantum error correction, researchers consider the
joint rank $(p,q)$-matricial range of $\bA$, denoted by
$\Lambda_{p,q}(\bA)$, consisting of
$m$-tuples of $k\times k$ matrices
$\bB = (B_1, \dots, B_m)$ such that for some unitary $U$,
$I_p\otimes B_j$ is the leading principal submatrix of
$U^*A_jU$ for $j = 1, \dots, m$; see \cite{LLPS} and its references.
One can define the $\Lambda_{p,q}^{ess}(\bA)$ as
$$\Lambda_{p,q}^{ess}(\bA)
= \cap \{\cl(\Lambda_{p,q}(\bA + \bK)): \bK \in \cS_\cK(\cH)^m\},$$
and  deduce the following results from those in the previous sections.

\newpage
\begin{theorem} \label{thm4.1} Let $\bA \in \cS(\cH)^m$.
\begin{itemize}
\item[{\rm (1)}] For any positive integer $N$, there is
$\bK\in \cS_\cK(\cH)^m$ such that
$$\Lambda_{p,q}(\bA + \bK) =
\Lambda^{ess}_{p,q}(\bA),$$
which is the same as
$$W_{ess}^q(\bA)
= W^q(\pi(A)) = W^q( \bA + \bK)
= W^q_s( \bA + \bK)
\quad \hbox{ for all } q \in \{1, \dots, N\}.$$
\item[{\rm (2)}]
Suppose $W_{ess}(\bA)$ is a simplex in $\IR^{1\times m}$. 
Then there is $\bK\in \cS_\cK(\cH)^m$
such that $$\Lambda_{p,q}(\bA+\bK) = \Lambda^{ess}_{p,q}(\bA),$$
which is the same as
$$W_{ess}^q(\bA)
=  W^q_s( \bA + \bK) =  W^q(\pi(A)) = W^q( \bA + \bK)$$
for all positive integer $q$.
\end{itemize}
\end{theorem}

\it Proof. \rm
By the result in \cite{LLPS},
$$\Lambda_{p,q}^{ess}(\bA) = W_{ess}^q(\bA) \subseteq
\Lambda_{p,q}(\bA + \bK) \subseteq W^q_s(\bA+\bK).$$
The results follow readily from Theorem \ref{thm3.1} and Theorem \ref{thm3.3}. \qed

\section*{Acknowledgment}

Li is an affiliate member of the Institute for Quantum Computing, University of
Waterloo, and an honorary professor of the Shanghai University; 
his research was supported by the USA NSF DMS 1331021, the Simons Foundation
Grant 351047, and NNSF of China Grant 11571220.
Paulsen is a professor of Pure Mathematics and the Institute for Quantum Computing, University of Waterloo, his research is partially supported by NSERC 03784. Poon would like to express his gratitude to  the generous support from   the Institute for Quantum Computing during his visit in September to October, 2017.

\noindent
{\bf Addresses}

\noindent(C.K. Li)
Department of Mathematics, College of William \& Mary,
Williamsburg, VA 23187, U.S.A.
(ckli@math.wm.edu)

\noindent
(V.I. Paulsen)
Institute for Quantum Computing and Department of Pure Mathematics,
University of Waterloo, Waterloo, Ontario, Canada
(vpaulsen@uwaterloo.ca)

\noindent
(Y.T. Poon)
Department of Mathematics, Iowa State University,
Ames, IA 50011, U.S.A.\newline
Center for Quantum Computing, Peng Cheng Laboratory, Shenzhen, 518055, China
(ytpoon@iastate.edu)

\end{document}